\documentclass[10pt]{amsart}
\usepackage{latexsym}
\usepackage{amsmath}
\usepackage{amssymb}
\usepackage{amscd}
\usepackage{array}
\usepackage{hhline}
\usepackage{color}
\usepackage[all]{xy}
\usepackage{enumerate}
\usepackage{graphicx}
\usepackage{lscape}

\theoremstyle{break}
\newtheorem{de}{Definition}[section]
\newtheorem{thm}[de]{Theorem}
\newtheorem{pro}[de]{Proposition}
\newtheorem{rem}[de]{Remark}
\newtheorem{lem}[de]{Lemma}
\newtheorem{cor}[de]{Corollary}

\newtheorem{ex}[de]{Example}

\def\Q{{\mathbb{Q}}}
\def\Z{{\mathbb{Z}}}
\def\c{{\mathcal{C}}}
\def\deg{{\mathrm{deg_{\Z}}}}
\def\pdeg{{\mathrm{deg}}}
\def\dim{{\mathrm{dim}}}

\def\ob{{\mathrm{Ob}}}
\def\mor{{\mathrm{Mor}}}
\def\M{{\mathrm{M}}}
\def\H{{\mathrm{H}}}

\def\id{{\mathrm{Id}}}
\def\MF{{\mathrm{MF}}}
\def\HMF{{\mathrm{HMF}}}

\title{Matrix factorizations and double line in $\mathfrak{sl}_n$ quantum link invariant}
\author{Yasuyoshi Yonezawa}
\address{Graduate School of Mathematics, Nagoya University\\ 464-8602 Furocho, Chikusaku, Nagoya, Japan }
\email{m03039e@math.nagoya-u.ac.jp}
\dedicatory{Dedicated to Professor Akihiro Tsuchiya on the occasion of his retirement}
\date{}

\begin{document}
\maketitle
\begin{abstract}
This article gives matrix factorizations for the trivalent diagrams and double line appearing in $\mathfrak{sl}_n$ quantum link invariant.
These matrix factorizations reconstruct Khovanov-Rozansky homology.
And we show that the Euler characteristic of the matrix factorization for a double loop equals the quantum dimension of the representation $\land^2 V$ of 
$U_q (\mathfrak{sl}_n)$ in Section \ref{sec3.3}.
\end{abstract}
\tableofcontents
%
%
%
%

\section{Introduction}
L. Kauffman introduced a graphical link invariant which is the normalized Jones polynomial \cite{Jo}\cite{Kau}.
It is well-known that the polynomial is derived from the fundamental representation of the quantum group $U_q(\mathfrak{sl}_2)$.
Further, G. Kuperberg constructed a graphical link invariant associated 
with the fundamental representation of the quantum group $U_q(\mathfrak{sl}_3)$  \cite{Kup}.
H. Murakami, T. Ohtsuki and S. Yamada introduced a graphical regular link invariant 
for the fundamental representation of the quantum group $U_q(\mathfrak{sl}_n)$ \cite{MOY} .
In general, we can also obtain a graphical link invariant for a given quantum group $U_q(\mathfrak{g})$ and the fundamental representation.
These invariants are collectively called $\mathfrak{g}$ quantum link invariants.\\
\indent
M. Khovanov constructed a categorification of $\mathfrak{sl}_2$ quantum link invariant \cite{K1}.
A categorification generally means the replacement of a set with a category by corresponding an element to an object. 
The morphism of the category is properly chosen to carry theory well-done. 
For a categorification, there is an inverse operation called a decategorification which is the replacement of a category with a set.
The decategorification of equivalent objects in the category is a same element in the set.\\
\indent
The Khovanov's theory is a beautiful example of the categorification;
this is the replacement of Jones polynomial $J(L)$, which is a map from the set of links to ``$\Z [q,q^{-1}]$'' ,
by a map $C_K$ from the set of links to ``the homotopy category of the bounded complex of graded $\Z$-modules''.
The bounded complex $C_K(L)$ and the $\Z \oplus \Z$-graded homology groups $KH^{i,j} (L)$ associated with $\c_K (L)$ 
also become link invariants under the Reidemeister moves. 
The decategorification is a $\Z$-graded Euler characteristic $\chi_{KH}$ with the normalized Jones polynomial $J(L)$ ;
$$
\chi_{KH} \left( \bigoplus_{i,j \in \Z} KH^{i,j} (L)\right) := \sum_{i,j \in \Z} (-1)^i \dim_{\Q} (KH^{i,j} (L) \otimes \Q ) q^j
= J(L).
$$
\indent
Recently, M. Khovanov and L. Rozansky introduced a categorification of $\mathfrak{sl}_n$ quantum link invariant 
using the homotopy category of the bounded complex of matrix factorizations \cite{KR1}\cite{KR2}. 
Since a resolution of a link diagram consists of a combination of the two local diagrams \input{figsmoothing4sln3} and \input{figsmoothing2sln} 
(See FIGURE \ref{fig-resolution}),
a matrix factorization for the resolution is defined by a tensor product of some matrix factorizations 
for the two local diagrams. 
Then, the categorification of $\mathfrak{sl}_n$ quantum link invariant is constructed 
as the map $\c_{KR}$ from the set of links to the homotopy category of the bounded complex of matrix factorizations for some resolutions. 

The bounded complex $\c_{KR}$ and the $\Z \oplus \Z \oplus \Z_2$-graded homology groups $KRH^{i,j,k} (L)$ associated with $\c_{KR} (L)$ 
also become link invariants under the Reidemeister moves.
The Euler characteristic $\chi_{KRH}$ is defined by
$$
\chi_{KRH} \left( \bigoplus_{i,j \in \Z , k \in \Z_2} KRH^{i,j,k} (L)\right) 
:= \sum_{i,j \in \Z , k \in \Z_2} (-1)^i \dim_{\Q} (KRH^{i,j,k} (L) ) q^j .
$$
This equals $\mathfrak{sl}_n$ quantum link invariant $\left< L \right>_n$ for the link $L$.
\begin{figure}[htb]
\begin{eqnarray*}
\Big\langle \ \input{figplus}\ \Big\rangle_{n}&=& \ q^{n-1} \Big\langle \input{figsmoothing1sln} \Big\rangle_{n}  -q^{n}\Big\langle \input{figsmoothing2sln1} \Big\rangle_{n} \label{slnplus}\\
\Big\langle \ \input{figminus}\ \Big\rangle_{n}&=&\ q^{1-n}\Big\langle \input{figsmoothing1sln} \Big\rangle_{n} -q^{-n}\Big\langle \input{figsmoothing2sln1} \Big\rangle_{n} \label{slnminus} \\
\Big\langle \ D \sqcup D\acute{} \ \Big\rangle_{n} &=& \Big\langle D \Big\rangle_{n} \ \Big\langle D\acute{} \Big\rangle_{n} \hspace{1cm} {\rm for\ any\ diagrams\ }D \ {\rm and} \ D\acute{}  \label{sln}\\
\Big\langle \ \input{figcircle1} \ \Big\rangle_{n} &=& [n] \label{loop}\\
\Big\langle \ \input{figsmoothbubble1}\ \Big\rangle_{n} &=& [2] \ \Big\langle \hspace{0.4cm} \input{figsmoothing3sln1} \hspace{0.4cm} \Big\rangle_{n} \label{babble1}\\
\Big\langle \ \input{figsmoothbubble2}\ \Big\rangle_{n} &=& [n-1] \Big\langle \hspace{0.2cm} \input{figsmoothing4sln1} \hspace{0.2cm} \Big\rangle_{n} \label{babble2}\\
\Big\langle \ \input{figsmoothbubble3}\ \Big\rangle_{n} &=& \ \Big\langle \ \input{figsmoothing4sln}\ \Big\rangle_{n} + [n-2]\ \Big\langle \input{figsmoothing3sln} \Big\rangle_{n} \label{babble3}\\
\Big\langle \ \input{figsmoothingreid1}\ \Big\rangle_{n} + \Big\langle \ \input{figsmoothingreid2}\ \Big\rangle_{n} &=& \Big\langle \ \input{figsmoothingreid3}\ \Big\rangle_{n} + \Big\langle \ \input{figsmoothingreid4}\ \Big\rangle_{n}  \label{babble4}
\end{eqnarray*}
where $[n]$ is a quantum integer $\frac{q^{n} - q^{-n}}{q-q^{-1}}$.
\caption{$\mathfrak{sl}_n$ quantum link invariant}
\label{fig-sln-skein}
\end{figure}
\\
\begin{figure}[htb]
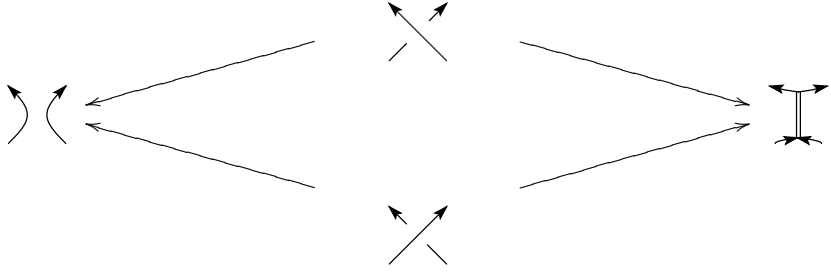

$$
\xymatrix{&&*^{\hspace{2cm}\input{figminus}\hspace{2cm}}\ar[rrd]\ar[lld] &&\\
*^{\input{figsmoothing1sln}\hspace{0.4cm}}&&*^{}&&*^{\hspace{0.4cm} \input{figsmoothing2sln1}}\\
&&*^{\hspace{2cm}\input{figplus}\hspace{2cm}}\ar[rru]\ar[llu] &&
}
$$
\caption{Two type resolutions for crossing}
\label{fig-resolution}
\end{figure}
\\
\indent
In this paper, we define matrix factorizations for the more general local diagrams 
\input{fig3valent-in} , \input{fig3valent-out} and the double line \input{figdline2}.
And we show that these matrix factorizations have some suitable properties; Proposition \ref{tetra} and corollary \ref{dloop} in Section \ref{sec3.3}.\\
\indent
Proposition \ref{tetra} claims that the inner marking made by gluing these local diagrams can be removed in the homotopy category of a matrix factorization.
For example, the matrix factorization for the diagram \input{figsmoothing2sln} 
is obtained as a tensor product of two matrix factorizations for \input{fig3valent-in} and \input{fig3valent-out}.\\
\indent
M. Khovanov and L. Rozansky showed that the Euler characteristic of a matrix factorization for the single loop 
\unitlength 0.1in
\begin{picture}(2.750,1.50)(2.2500,-4.0)
%
\special{pn 8}%
\special{ar 350 350 120 120  0.0000000 6.2831853}%
\end{picture}%
 equals to 
the quantum dimension of the fundamental representation $V$ of $U_q(\mathfrak{sl}_n)$.
Corollary \ref{dloop} claims that the Euler characteristic of a matrix factorization for the double loop 
\unitlength 0.1in
\begin{picture}(3.30,1.500)(1.85,-3.75)
%
\special{pn 8}%
\special{ar 350 350 120 120  0.0000000 6.2831853}%
%
\special{pn 8}%
\special{ar 350 350 100 100  0.0000000 6.2831853}%
\end{picture}%
 also equals to 
the quantum dimension of $\land^2 V$ of $U_q(\mathfrak{sl}_n)$.\\
\\ \\
{\bf Acknowledgements:} 
The author would like to thank Hiroaki Kanno, Hiroyuki Ochiai, Hidefumi Ohsugi and Akihiro Tsuchiya for many helpful discussion. 
And he would also like to thank 
Mikhail Khovanov and Lev Rozansky for supplying their paper \cite{KR3} (preliminary version) to him. 
%
%
%
%

\section{Category of matrix factorization}
\subsection{Matrix factorization}
We describe a category of a matrix factorization. M. Khovanov and L. Rozansky first imported this algebraic object into link theory \cite{KR1}\cite{KR2}.\\
\indent
Let $R$ be a polynomial ring over $\Q$ and let $M_0$,$M_1$ be free $R$-modules permitted infinite rank. $\overline{M}=(M_0,M_1)$ is a {\it matrix factorization} 
with a potential $\omega \in R$ if $\overline{M}$ consists of a $2$-cyclic complex of free $R$-modules $M_0$,$M_1$;
$$
\Big( \xymatrix{
*^{M_0}\ar[rrr]^{d_{M_0}}&&&*^{M_1}\ar[rrr]^{d_{M_1}}&&&*^{M_0}
}
\Big)
$$
such that $d_{M_1} d_{M_0} = \omega\,\id_{M_0}$ and $d_{M_0} d_{M_1} = \omega\,\id_{M_1}$.\\
\indent
For a polynomial ring $R$ and a polynomial $\omega \in R$, let $\MF_{R,\omega}$ be a category of a matrix factorization 
whose object is a matrix factorization $(M_0,M_1)$ with the potential $\omega$ 
and whose morphism $\overline{f}=(f_0,f_1)$ between $\overline{M}=(M_0,M_1)$ and $\overline{N}=(N_0,N_1)$ consists of a pair of $R$-module morphisms 
$f_0:M_0 \to N_0$ and $f_1:M_1 \to N_1$ such that $d_{N_0} f_0 = f_1 d_{M_0}$ and 
$d_{N_1} f_1 = f_0 d_{M_1}$. \\
%
\indent
For $\overline{M} \in \ob(\MF_{R,\omega})$ and $\overline{N} \in \ob(\MF_{R\acute{},\omega\acute{}})$, 
we define the tensor product $\overline{M}\boxtimes \overline{N} \in \ob(\MF_{R\otimes R\acute{},\omega + \omega \acute{}})$ by
$$
\overline{M}\boxtimes \overline{N} :=
\left(
\xymatrix{
*^{\left(
	\begin{array}{c}
		M_0\otimes N_0\\
		M_1\otimes N_1
	\end{array}
    \right)}
\ar[rrr]_{\left(
		\begin{array}{cc}
			d_{M_0}&-d_{N_1}\\
			d_{N_0}&d_{M_1}
		\end{array}
	    \right)}&&&
*^{\left(
	\begin{array}{c}
		M_1\otimes N_0\\
		M_0\otimes N_1
	\end{array}
    \right)}
\ar[rrr]_{\left(
		\begin{array}{cc}
			d_{M_1}&d_{N_1}\\
			-d_{N_0}&d_{M_0}
		\end{array}
	    \right)}
&&&*^{\left(
	\begin{array}{c}
		M_0\otimes N_0\\
		M_1\otimes N_1
	\end{array}
    \right)}
}
\right)
,
$$
where 
\begin{center}
$\left(
		\begin{array}{cc}
			d_{M_0}&-d_{N_1}\\
			d_{N_0}&d_{M_1}
		\end{array}
	    \right)$ 
and
$\left(
		\begin{array}{cc}
			d_{M_1}&d_{N_1}\\
			-d_{N_0}&d_{M_0}
		\end{array}
	    \right) $
\end{center}
simply denote 
\begin{center}
$\left(
		\begin{array}{cc}
			d_{M_0} \otimes \id_{N_0}&- \id_{M_1}\otimes d_{N_1}\\
			\id_{M_0}\otimes d_{N_0}&d_{M_1}\otimes \id_{N_1}
		\end{array}
	    \right)$ 
and
$\left(
		\begin{array}{cc}
			d_{M_1}\otimes \id_{N_0}&\id_{M_0}\otimes d_{N_1}\\
			-\id_{M_1}\otimes d_{N_0}&d_{M_0}\otimes \id_{N_1}
		\end{array}
	    \right) $.
\end{center}
\begin{rem}
We consider $R \otimes R\acute{}$ as a tensor product of polynomial rings $R$ and $R\acute{}$ over $R \cap R\acute{}$.
We also consider $M \otimes M\acute{}$ as a tensor product of an $R$-module $M$ and an $R\acute{}$-module $M\acute{}$ over $R \cap R\acute{}$.
\end{rem}
\begin{lem}\label{com-ass}
{\rm (1)}For $\overline{M} \in \ob(\MF_{R,\omega})$ and $\overline{N} \in \ob(\MF_{R\acute{},\omega\acute{}})$, 
there is an isomorphism in $\MF_{R\otimes R\acute{},\omega + \omega \acute{}}$
$$\overline{M}\boxtimes \overline{N} \simeq \overline{N}\boxtimes \overline{M}.$$ 
{\rm (2)}For $\overline{L} \in \ob(\MF_{R,\omega})$, $\overline{M} \in \ob(\MF_{R\acute{},\omega\acute{}})$ 
and $\overline{N} \in \ob(\MF_{R \acute{}\,\acute{},\omega \acute{}\,\acute{}})$, 
there is an isomorphism in $\MF_{R\otimes R\acute{}\otimes R \acute{}\,\acute{},\omega + \omega\acute{} +\omega \acute{}\,\acute{}}$ 
$$(\overline{L}\boxtimes \overline{M})\boxtimes \overline{N} \simeq \overline{L}\boxtimes (\overline{M}\boxtimes \overline{N}).$$ 
\end{lem}
\begin{proof}
(1)For $R$-modules $M$ and $N$, we define $T:M\otimes N \to N \otimes M$ by $T(m \otimes n) := n \otimes m$.
$
\overline{f}=
\left(
\left(
	\begin{array}{cc}
		T&0\\
		0&-T
	\end{array}
\right),
\left(
	\begin{array}{cc}
		0&T\\
		T&0
	\end{array}
\right)
\right)
$
is a morphism from $\overline{M}\boxtimes \overline{N}$ to $\overline{N}\boxtimes \overline{M}$ 
and also a morphism from $\overline{N}\boxtimes \overline{M}$ to $\overline{M}\boxtimes \overline{N}$.
Since $\overline{f}^2 = \overline{\id}$, $\overline{f}$ gives isomorphic between these matrix factorizations.\\
(2)By definition, we have
$$d_{((\overline{L}\boxtimes \overline{M})\boxtimes \overline{N})_0}=
\left(
	\begin{array}{cccc}
		d_{L_0}&-d_{M_1}&-d_{N_1}&0\\
		d_{M_0}& d_{L_1}&0       &-d_{N_1}\\
		d_{N_0}&0       & d_{L_1}&d_{M_1}\\
		0      & d_{N_0}&-d_{M_0}&d_{L_0}
	\end{array}
\right)
,
d_{((\overline{L}\boxtimes \overline{M})\boxtimes \overline{N})_1}=
\left(
	\begin{array}{cccc}
		 d_{L_1}& d_{M_1}&d_{N_1}&0\\
		-d_{M_0}& d_{L_0}&0       &d_{N_1}\\
		-d_{N_0}&0       & d_{L_0}&-d_{M_1}\\
		 0      &-d_{N_0}& d_{M_0}&d_{L_1}
	\end{array}
\right)
$$ and
$$d_{(\overline{L}\boxtimes (\overline{M}\boxtimes \overline{N}))_0}=
\left(
	\begin{array}{cccc}
		d_{L_0}&0       &-d_{M_1}&-d_{N_1}\\
		0      & d_{L_0}& d_{N_0}&-d_{M_0}\\
		d_{M_0}&-d_{N_1}& d_{L_1}&0\\
		d_{N_0}& d_{M_1}&0       &d_{L_1}
	\end{array}
\right)
,
d_{(\overline{L}\boxtimes (\overline{M}\boxtimes \overline{N}))_1}=
\left(
	\begin{array}{cccc}
		 d_{L_1}&0       & d_{M_1}&d_{N_1}\\
		 0      & d_{L_1}&-d_{N_0}&d_{M_0}\\
		-d_{M_0}& d_{N_1}& d_{L_0}&0\\
		-d_{N_0}&-d_{M_1}&0       &d_{L_0}
	\end{array}
\right).
$$\\
Thus it is obvious that $\overline{f}=\left(
\left(
	\begin{array}{cccc}
		 1&0&0&0\\
		 0&0&0&1\\
		 0&1&0&0\\
		 0&0&1&0
	\end{array}
\right)
,
\left(
	\begin{array}{cccc}
		 1&0&0&0\\
		 0&0&0&1\\
		 0&1&0&0\\
		 0&0&1&0
	\end{array}
\right)
\right) :(\overline{L}\boxtimes \overline{M})\boxtimes \overline{N} \to \overline{L}\boxtimes (\overline{M}\boxtimes \overline{N})$ 
is an isomorphism.
\end{proof}
\begin{lem}\label{identity}
The matrix factorization $(\xymatrix{R \ar[r]^0&0\ar[r]^0&R})$ is the unit object in $\MF_{R,\omega}$. That is, 
for any matrix factorization $\overline{M} \in \MF_{R,\omega}$,
$$
\overline{M} \boxtimes (\xymatrix{R \ar[r]^0&0\ar[r]^0&R}) \simeq \overline{M} .
$$
\end{lem}
\begin{proof}
By definition, we have
\begin{eqnarray*}
\overline{M}\boxtimes (\xymatrix{R \ar[r]^{0}&0\ar[r]^{0}&R})&=&
\left(
\xymatrix{
*^{\left(
\begin{array}{c}
M_0 \otimes R\\
M_1 \otimes 0
\end{array}
\right) }
\ar[rr]^{
\left(
\begin{array}{cc}
d_{M_0}&0\\
0&d_{M_1}
\end{array}
\right) }&&
*^{\left(
\begin{array}{c}
M_1 \otimes R\\
M_0 \otimes 0
\end{array}
\right)}
\ar[rr]^{
\left(
\begin{array}{cc}
d_{M_1}&0\\
0&d_{M_0}
\end{array}
\right) }&&
*^{\left(
\begin{array}{c}
M_0 \otimes R\\
M_1 \otimes 0
\end{array}
\right)}}
\right)
\\
&\simeq&\overline{M}.
\end{eqnarray*}
\end{proof}

\indent
The translation functor $\left< 1 \right>$ changes the matrix factorization $\overline{M} = (M_{0},M_{1})$ into 
$$
\overline{M} \left< 1 \right> =
\Big(
\xymatrix{
*{M_{1}}\ar[rrr]^{-d_{M_1}}&&&*{M_{0}}\ar[rrr]^{-d_{M_0}}&&&*{M_{1}}
}
\Big)
.
$$
The functor $\left< 2 \right> (= \left< 1 \right>^2)$ is the identity functor.
\begin{lem}
For $\overline{M} \in \ob(\MF_{R,\omega})$ and $\overline{N} \in \ob(\MF_{R\acute{},\omega\acute{}})$, 
there is an isomorphism in $\MF_{R\otimes R\acute{},\omega + \omega \acute{}}$
\begin{eqnarray*}
(\overline{M}\boxtimes \overline{N})\left< 1 \right> &=& (\overline{M}\left< 1 \right> )\boxtimes \overline{N}\\
&\simeq& \overline{M}\boxtimes (\overline{N}\left< 1 \right> ).
\end{eqnarray*}
\end{lem}
\begin{proof}
We directly find that $(\overline{M}\boxtimes \overline{N})\left< 1 \right>$ equals $(\overline{M}\left< 1 \right> )\boxtimes \overline{N}$ 
by definition. The second equivalence is correct by Lemma \ref{com-ass} (1) and the first equality.
\end{proof}
\indent
The morphism $\overline{f}: \overline{M} \to \overline{N} \in \mor (\MF_{R,\omega})$ is {\it null-homotopic} 
if morphisms $h_{0}:M_0 \to N_1$ and $h_{1}:M_1 \to N_0$ exist 
such that $f_0 = h_1 d_{M_0} + d_{N_1} h_0$ and $f_1 = h_0 d_{M_1} + d_{N_0} h_1$.
And $\overline{f},\overline{g}: \overline{M} \to \overline{N} \in \mor (\MF_{R,\omega})$ are {\it homotopic} 
if $\overline{f} - \overline{g}$ is null-homotopic.\\
\indent
Let $\HMF_{R,\omega}$ be the quotient category of $\MF_{R,\omega}$ which has the same objects to $\ob(\MF_{R,\omega})$
and has morphisms of $\mor(\MF_{R,\omega})$ modulo null-homotopic. A matrix factorization in $\MF_{R,\omega}$ is called {\it contractible} 
if it is isomorphic to the zero matrix factorization
$$
\Big(
\xymatrix{
*{0}\ar[rrr]^{0}&&&*{0}\ar[rrr]^{0}&&&*{0} 
}
\Big)
$$
in $\HMF_{R,\omega}$. 
\begin{ex}
Let $R_0$,$R_1$ be a ring $R$ and $\omega \in R$.
$$
\Big(
\xymatrix{
*{R_{1}}\ar[rrr]^{\omega}&&&*{R_{0}}\ar[rrr]^{1}&&&*{R_{1}} 
}
\Big)
$$
and
$$
\Big(
\xymatrix{
*{R_{1}}\ar[rrr]^{1}&&&*{R_{0}}\ar[rrr]^{\omega}&&&*{R_{1}} 
}
\Big)
$$
are contractible.
\end{ex}
%
%
%
%
\subsection{$\Z$-graded matrix factorization}\label{sec2.2}
Let $R$ be replaced with a $\Z$-graded polynomial ring over $\Q$ whose each parameter has a $\Z$-grading 
and let $M_0$, $M_1$ be also replaced with free $\Z$-graded $R$-modules.
The category $\MF^{gr}_{R,\omega}$ is the category of a $\Z$-graded matrix factorization whose object is 
the same object to $\ob (\MF_{R,\omega})$ except having a $\Z$-grading 
and whose morphism consists of a morphism with preserving a $\Z$-grading in $\mor (\MF_{R,\omega})$. 
The quotient category $\HMF^{gr}_{R,\omega}$ is also defined 
as $\ob (\HMF^{gr}_{R,\omega}) = \ob (\MF^{gr}_{R,\omega})$ and $\mor (\HMF^{gr}_{R,\omega}) = \mor (\MF^{gr}_{R,\omega}) / $\{ null-homotopic \} . \\
\indent
The $\Z$-grading shift $\{ n \}$ ($n \in \Z$) turns the matrix factorization $\overline{M}$ into 
$$
\Big( \xymatrix{ M_0\{ n \} \ar[r]^{ d_{M_0} }   &   M_1\{ n\}   \ar[r]^{ d_{M_1} }   &   M_0\{ n\} }\Big)
.
$$
\begin{lem}
For $\overline{M} \in \ob(\MF^{gr}_{R,\omega})$ and $\overline{N} \in \ob(\MF^{gr}_{R\acute{},\omega\acute{}})$, 
there is an equality in $\MF^{gr}_{R\otimes R\acute{},\omega + \omega \acute{}}$
\begin{eqnarray*}
(\overline{M}\boxtimes \overline{N})\{ n\} &=& (\overline{M}\{ n\} )\boxtimes \overline{N}\\
&=& \overline{M}\boxtimes (\overline{N}\{ n\} )
\end{eqnarray*}
\end{lem}
\begin{proof}
We find that these objects are really identical by definition. 
\end{proof}
\indent
For the $\Z$-graded matrix factorization $\overline{M}$ in $\MF^{gr}_{R,0}$, 
we can consider $\Z \oplus \Z_2$-graded homology group as follows;
$$ \H(\overline{M})=\bigoplus_{j \in \Z ,k\in \Z_2 } \H^{j,k}(\overline{M})$$
whose $k$ is a complex grading of the matrix factorization $\overline{M}$, i.e. $i=0$ or $1$, 
and $j$ is a $\Z$-grading induced by the $\Z$-graded modules of the matrix factorization $\overline{M}$.
The {\it Euler characteristic} $\overline{\chi}$ is defined by
\begin{equation*}
\overline{\chi}(\H(\overline{M})) = \sum_{ j \in \Z, k \in \Z_2} \dim_{\Q} \H^{j,k}(\overline{M})  q^{j} \label{euler}.
\end{equation*}

%
%
%
%

\subsection{Koszul matrix factorization}
Let $R$ be a $\Z$-graded polynomial ring over $\Q$ and let $\deg (a)$ denote a $\Z$-grading of the polynomial $a \in R$. 
For polynomials $a$, $b \in R$ and a $\Z$-graded $R$-module $M$, 
we define the matrix factorization $ K(a;b)_{M} $
with the potential $a b$ by
$$
K(a;b)_{M}:= 
\Big(
\xymatrix{
*{M}\ar[rrr]^(.3){a}&&&*{M\{ \frac{1}{2}(\deg (b)-\deg (a))\}}\ar[rrr]^(.7){b}&&&*{M}
}
\Big)
.$$ 
\begin{rem}
Let $R$ be a polynomial ring over $\Q$ and let $R_y$ be the polynomial ring $R[y]$. For polynomials $a$ and $b$ in $R$,
$K(a;b)_{R_y}$ is a matrix factorization of $R_y$-modules with rank $1$ as an object in $\MF^{gr}_{R_y,ab}$.
And we can consider that $K(a;b)_{R_y}$ is a matrix factorization of $R$-modules with infinite rank as an object in $\MF^{gr}_{R,ab}$.
\end{rem}
\begin{lem}\label{functor}
$$K(a;b)_{M}\,\left< 1 \right> = K(-b;-a)_{M}\,\{ \frac{1}{2}(\deg b -\deg a) \} .$$
\end{lem}
\begin{proof}
By definition, we have 
\begin{eqnarray*}
K(a;b)_{M}\left< 1 \right> &=& \Big(\xymatrix{M \ar[rr]^(.3){a}&& M\{ \frac{1}{2}(\deg b -\deg a) \} \ar[rr]^(.7){b}&& M}\Big)\,\left< 1 \right> \\
&=& \Big(\xymatrix{M\{ \frac{1}{2}(\deg b -\deg a) \} \ar[rr]^(.7){-b}&& M \ar[rr]^(.3){-a}&& M\{ \frac{1}{2}(\deg b -\deg a) \} }\Big)\\
&=& \Big(\xymatrix{M \ar[rr]^(.3){-b}&& M\{ \frac{1}{2}(\deg a -\deg b) \} \ar[rr]^(.7){-a}&& M}\Big)\,\{ \frac{1}{2}(\deg b -\deg a) \} \\
&=& K(-b;-a)_{M}\,\{ \frac{1}{2}(\deg b -\deg a) \} .
\end{eqnarray*}
\end{proof}
In general, for the sequences $\mathbf{a}={}^t(a_1, a_2, \ldots, a_k)$ and $\mathbf{b}={}^t(b_1, b_2, \ldots, b_k)$ of polynomials in $R$, 
we define the matrix factorization $K\left( \mathbf{a} ; \mathbf{b} \right)_{M}$ with the potential $\sum_{i=1}^k a_i b_i$ by 
\[
K\left( \mathbf{a} ; \mathbf{b} \right)_{M}
:=
\mathop{\boxtimes}_{i=1}^k K(a_i;b_i)_{M}
.
\]
This matrix factorization is called a {\it Koszul matrix factorization}.
\begin{lem} \label{equiv}
Let 
$c$ be a non-zero element in $\Q$.\\

There is an isomorphism in $\MF^{gr}_{R,a b}$
$$
K(a;b)_M \simeq K(c a;c^{-1} b)_M. 
$$

\end{lem}

\begin{proof}
$\overline{f}=(1,c):K(a;b)_M \to K(c a ; c^{-1} b)_M$ and $\overline{g}=(1,c^{-1}):K(c a ; c^{-1} b)_M \to K(a;b)_M$ satisfy
 that $\overline{g} \, \overline{f} = \overline{\id}$ and $\overline{f} \, \overline{g} = \overline{\id}$.
Thus, these matrix factorizations are isomorphic.\\
\end{proof}

\begin{thm}\label{exclude}{\rm [Khovanov-Rozansky,Theorem 2.2.\cite{KR3}]}
We put $R=\Q [x_1,x_2,\cdots ,x_n]$ and $R_y = R [y]$. Let ${}^t( a_1 , a_2 , \cdots , a_k )$ and 
${}^t( b_1 , b_2 , \cdots , b_k )$ be sequences of polynomials $a_j , b_j \in R_y$.\\
We assume that $K(\mathbf{a};\mathbf{b})_{R_y}$ 
is a Koszul matrix factorization with the potential $\omega \in R$. 
That is to say, we can seem that this matrix factorization is an object in $HMF^{gr}_{R,\omega }$, which consists of infinite rank $R$-modules $R_y$.
Furthermore, we assume that $b_i =c \, y^n + p$ for some $i$, where $c$ is a non-zero element in $\Q$ 
and p is the polynomial in $R_y$ whose degree for $y$ is less than n. Then, there is an isomorphism in $\HMF^{gr}_{R,\omega }$
$$K(\mathbf{a};\mathbf{b})_{R_y} \simeq K(\stackrel{i}{\check{\mathbf{a}}};\stackrel{i}{\check{\mathbf{b}}})_{R_y/\left< b_i\right> },$$
where $\stackrel{i}{\check{\mathbf{a}}}$ and $\stackrel{i}{\check{\mathbf{b}}}$ are associated with 
$\mathbf{a}$ and $\mathbf{b}$ removing the $i$-th polynomial.
\end{thm}

\begin{proof}
We can replace $b_i$ with $y^n + p$ using Lemma \ref{equiv} .
Then, we repeat proof by M.Khovanov and L.Rozansky in \cite{KR3}.
By definition and the above assumption, we have 
$$
K(\mathbf{a};\mathbf{b})_{R_y} = K(\stackrel{i}{\check{\mathbf{a}}};\stackrel{i}{\check{\mathbf{b}}})_{R_y} \boxtimes K(a_i;y^n+p)_{R_y}.
$$
The Koszul matrix factorization $K(\stackrel{i}{\check{\mathbf{a}}};\stackrel{i}{\check{\mathbf{b}}})_{R_y}$ is described as
$$
\Big(
\xymatrix{
R_{y}^r \ar[rrr]^{D_0} &&& R_{y}^r \ar[rrr]^{D_1} &&& R_{y}^r,
}
\Big) ,
$$
where $D_0$, $D_1 \in \M_r(R_y)$ and $r=2^{k-2}$.
Thus, we obtain
$$
K(\mathbf{a};\mathbf{b})_{R_y}=
\left(
\xymatrix{
*{
\begin{array}{c}
R_y^r\\
\oplus\\
R_y^r
\end{array}
}
\ar[rrr]^{D_0 \acute{}}&&&
*{
\begin{array}{c}
R_y^r\\
\oplus\\
R_y^r
\end{array}
}
\ar[rrr]^{D_1 \acute{}
}&&&
*{
\begin{array}{c}
R_y^r\\
\oplus\\
R_y^r
\end{array}
}
}
\right)
,$$
where $D_0 \acute{} =\left(
\begin{array}{cc}
D_0&(-y^n-p) \id_{R_y^r}\\
a_i \id_{R_y^r}&D_1
\end{array}
\right)$
,
$D_1 \acute{}=\left(
\begin{array}{cc}
D_1&(y^n+p) \id_{R_y^r}\\
-a_i \id_{R_y^r}&D_0
\end{array}
\right)$.

\indent
The ring $R_y$ is split into the direct sum as an $R$-module as follows;
\begin{center}
$R_y \simeq R_{<n} \oplus R_{\geq n},$ 
where $\displaystyle R_{<n}= \bigoplus_{i=0}^{n-1} R\cdot y^i $ 
and $\displaystyle R_{\geq n}=\bigoplus_{i=0}^{\infty} R \cdot y^i (y^n+p) .$
\end{center}
The $R$-module morphism $R_y \stackrel{y^n + p}{\longrightarrow} R_y$ induces the $R$-module isomorphism
$$
\xymatrix{
f_{iso}:R_y \ar[rr]^{y^n+p} && R_{\geq n}.
}
$$
Moreover, there are the natural $R$-module injection
$$
\xymatrix{
f_{inj}:R_{<n} \ar[rr] && R_y,
}
$$
and the natural $R$-module projections
$$
\xymatrix{
f_{proj_{<n}}:R_y \ar[rr] && R_{<n}
}
$$
and
$$
\xymatrix{
f_{proj_{\geq n}}:R_y \ar[rr] && R_{\geq n}.
}
$$

\indent
The $R$-module $R_y^r$ is also split into the direct sum 
$$
R_y^r \simeq R_{<n}^r \oplus R_{\geq n}^r.
$$
Then, there are the $R$-module isomorphism
$$
\xymatrix{
F_{iso}:R_y^r \ar[rr]^{(y^n+p)\id_{R^r_y}} && R_{\geq n}^r,
}
$$
the $R$-module injection
$$
\xymatrix{
F_{inj}:R_{<n}^r \ar[rr] && R_y^r
}
$$
and the $R$-module projections
$$
\xymatrix{
F_{proj_{<n}}:R_y^r \ar[rr] && R_{<n}^r
}
$$
and
$$
\xymatrix{
F_{proj_{\geq n}}:R_y^r \ar[rr] && R_{\geq n}^r .
}
$$

It is easy to find that the $R$-module morphisms
$$
\phi_0 =\xymatrix{
*{\begin{array}{c}
R_{<n}^r\\
\oplus\\
R_{\geq n}^r\\
\oplus\\
R_y^r
\end{array}}
\ar[rrrrrr]_{\left(
\begin{array}{ccc}
F_{inj}&(y^n+p)F_{iso}^{-1}&0\\
0&D_0 \, F_{iso}^{-1}&\id_{R_y^r}
\end{array}\right)
}&&&&&&
*{
\begin{array}{c}
R_y^r\\
\oplus\\
R_y^r
\end{array}
}
}$$
and 
$$
\phi_1 =\xymatrix{
*{\begin{array}{c}
R_{<n}^r\\
\oplus\\
R_{\geq n}^r\\
\oplus\\
R_y^r
\end{array}}
\ar[rrrrrr]_{\left(
\begin{array}{ccc}
F_{inj}&(-y^n-p)F_{iso}^{-1}&0\\
0&D_1 \, F_{iso}^{-1}&\id_{R_y^r}
\end{array}\right)
}&&&&&&
*{
\begin{array}{c}
R_y^r\\
\oplus\\
R_y^r
\end{array}
}
}$$
are $R$-module isomorphisms.
Let $\overline{M}$ be the matrix factorization
$$
\overline{M}=
\left(
\xymatrix{
*{\begin{array}{c}
R_{<n}^r\\
\oplus\\
R_{\geq n}^r\\
\oplus\\
R_y^r
\end{array}}
\ar[rrr]_{\phi_1^{-1}\, D_0 \acute{} \, \phi_0}&&&
*{\begin{array}{c}
R_{<n}^r\\
\oplus\\
R_{\geq n}^r\\
\oplus\\
R_y^r
\end{array}}\ar[rrr]_{\phi_0^{-1}\, D_1 \acute{} \, \phi_1}&&&
*{\begin{array}{c}
R_{<n}^r\\
\oplus\\
R_{\geq n}^r\\
\oplus\\
R_y^r
\end{array}}
}
\right)
.
$$
Then, $\overline{\phi}=$($\phi_0$,$\phi_1$) is an isomorphism from $\overline{M}$ to $K(\mathbf{a};\mathbf{b})_{R_y}$ in $\MF^{gr}_{R,\omega }$.\\
Since $\overline{\phi}^{-1}$ consists of 
$$\phi_0^{-1} = \left(
\begin{array}{cc}
F_{proj_{<n}}&0\\
(y^n+p) F_{iso}^{-1}\, F_{proj_{\geq n}}&0\\
-D_0 F_{iso}^{-1}\, F_{proj_{\geq n}}&\id_{R_y}^r
\end{array}
\right) {\rm and \,\,} \phi_1^{-1} = \left(
\begin{array}{cc}
F_{proj_{<n}}&0\\
(-y^n-p) F_{iso}^{-1}\, F_{proj_{\geq n}}&0\\
D_1 F_{iso}^{-1}\, F_{proj_{\geq n}}&\id_{R_y}^r
\end{array}
\right) ,$$
the morphisms $\phi_1^{-1}\, D_0 \acute{} \, \phi_0$ and $\phi_0^{-1}\, D_1 \acute{} \, \phi_1$ are described by
$$
\phi_1^{-1}\, D_0 \acute{} \, \phi_0=\left(
\begin{array}{ccc}
F_{proj_{<n}}\, D_0\, F_{inj} & 0 &0\\
\ast & 0 &F_{iso}\\
\ast & \omega F_{iso}^{-1} &0
\end{array}
\right) {\rm and} \,\,
\phi_0^{-1}\, D_1 \acute{} \, \phi_1=\left(
\begin{array}{ccc}
F_{proj_{<n}}\, D_1\, F_{inj} & 0 &0\\
\ast & 0 &F_{iso}\\
\ast & \omega F_{iso}^{-1} &0
\end{array}
\right).
$$
\indent
The matrix factorization obtained by restricting $R_{< n}^r \oplus R_{\geq n}^r \oplus R_y^r$ of $\overline{M}$ to $R_{\geq n}^r \oplus R_y^r$
$$
\left(
\xymatrix{
*^{
\left(
\begin{array}{c}
R_{\geq n}^r \\
\oplus \\
R_y^r
\end{array}
\right)
}
\ar[rrr]_{
\left(
\begin{array}{cc}
0 & F_{iso}\\
\omega F_{iso}^{-1} & 0
\end{array}
\right)
}
&&&
*^{
\left(
\begin{array}{c}
R_{\geq n}^r \\
\oplus \\
R_y^r
\end{array}
\right)
}
\ar[rrr]_{
\left(
\begin{array}{cc}
0 & F_{iso}\\
\omega F_{iso}^{-1} & 0
\end{array}
\right)
} &&&
*^{
\left(
\begin{array}{c}
R_{\geq n}^r \\
\oplus \\
R_y^r
\end{array}
\right)
}}
\right)
$$
is contractible in $\HMF^{gr}_{R,\omega }$.
Thus, $\overline{M}$ is isomorphic to the quotient matrix factorization
$$
\left(
\xymatrix{
*{R_{<n}^r}\ar[rrr]^{F_{proj_{<n}}\, D_0 \, F_{inj}}&&&*{R_{<n}^r}\ar[rrr]^{F_{proj_{<n}}\, D_1 \, F_{inj}}&&&*{R_{<n}^r}.
}
\right)
$$ in $\HMF^{gr}_{R,\omega }$. 
By the choice of a basis of $R_y$ as an $R$-module, it is easy to find that this quotient matrix factorization equals
 $K(\stackrel{i}{\check{\mathbf{a}}};\stackrel{i}{\check{\mathbf{b}}})_{R_y/\left< b_i \right> }$.
\end{proof}
%
%
%
%

\section{Matrix factorization for trivalent diagrams and double line}
\subsection{Khovanov Rozansky homology}

We briefly recall Khovanov-Rozansky link homology theory, only the definition of matrix factorizations for planar diagrams and 
some properties of the matrix factorizations.
See Section 6 in \cite{KR1} for further details.\\
\indent
In \cite{KR1}, M. Khovanov and L. Rozansky defined a categorification of planar graphs 
\,\,\input{figsmoothing2sln}\,\, and \,\,\,\input{figsmoothing4sln3}\,\,\, 
as a Koszul matrix factorization. \\
\indent
We assign a parameter $x_i$ on the end point of a single line. And each $\Z$-grading of parameters $x_i$ is $2$. 
This $\Z$-grading induces a $\Z$-grading of a matrix factorization in $\MF^{gr}_{R,\omega }$ and $\HMF^{gr}_{R,\omega }$. 
The function $f(s_1,s_2)$ is obtained by expanding the power sum $x^{n+1}+y^{n+1}$ with the elementary symmetric polynomials $x+y$ and $xy$, i.e.
the function $f(s_1,s_2)$ satisfies that
$$
f(x+y,xy)=x^{n+1}+y^{n+1}.
$$
\indent
We define matrix factorizations for the planar diagrams \input{figsmoothing4sln1-mf} and \input{figsmoothing2sln-mf} by
\begin{eqnarray*}
\c\acute{} \Bigg( \hspace{0.35cm} \input{figsmoothing4sln1-mf}\hspace{0.35cm} \Bigg)_{n}&=&K(\pi_{12};x_1-x_2)_{\Q [x_1,x_2]} ,\\
\c\acute{} \Bigg( \ \input{figsmoothing2sln-mf}\ \Bigg)_{n}&=&K \left(
\left(
\begin{array}{c}
u^{1,2}_{3,4}\\
v^{1,2}_{3,4}
\end{array}
\right) ;
\left(
\begin{array}{c}
x_1 + x_2 - x_3 - x_4\\
x_1 x_2 - x_3 x_4
\end{array}
\right)
\right)_{\Q [x_1,x_2,x_3,x_4]}\{ -1 \},
\end{eqnarray*}
where $\pi_{12} (x_1 - x_2) = x_1^{n+1}-x_2^{n+1}$, $u^{1,2}_{3,4}= \frac{f(x_1 + x_2,x_1 x_2)-f(x_3 + x_4,x_1 x_2)}{x_1 + x_2 - x_3 - x_4}$ 
and $v^{1,2}_{3,4}=\frac{f(x_3 + x_4,x_1 x_2)-f(x_3 + x_4,x_3 x_4)}{x_1 x_2 - x_3 x_4}$.
\indent
Furthermore, we construct a matrix factorization for the more complex planar diagrams 
produced by combinatorially joining the diagrams \,\,\input{figsmoothing2sln}\,\, and \,\,\,\input{figsmoothing4sln3}. 
Now, suppose the planar diagrams \input{figgluing-mf1} and \input{figgluing-mf2}, 
which can match at $x$ and $x\acute{}$ with keeping the orientation.
For the matrix factorizations 
$$
\c\acute{}\Big(\input{figgluing-mf1}\Big)_{n} \in \ob (\MF^{gr}_{R,\omega + x^{n+1}})
\hspace{0.5cm} {\rm and} \hspace{0.5cm} 
\c\acute{}\Big(\input{figgluing-mf2}\Big)_{n}
\in \ob (\MF^{gr}_{R\acute{},\omega\acute{} - x\acute{}^{n+1}}),
$$
 we define the matrix factorization 
$$
\c\acute{}\Bigg( \input{figgluing2} \Bigg)_{n} \in \ob (\MF^{gr}_{R \otimes R\acute{},\omega + \omega\acute{}}) 
\hspace{0.5cm} {\rm by} \hspace{0.5cm} 
\c\acute{}\Big(\input{figgluing-mf1}\Big)_{n} \boxtimes \c\acute{}\Big(\input{figgluing-mf2}\Big)_{n}\Big|_{x\acute{}=x}.
$$
This means that we identify the top parameter $x$ and the tail parameter $x\acute{}$ after taking the tensor product of these matrix factorizations.
And if a planar diagram consists of two disjoint planar diagrams $D_1$ and $D_2$,
then this matrix factorization is defined by the tensor product of matrix factorizations for $D_1$ and $D_2$;
$$
\c\acute (D_1 \bigsqcup D_2)_n = \c\acute (D_1)_{n} \boxtimes \c\acute (D_2)_{n} .
$$
\begin{pro}
The map $\c\acute{}$ from a planar diagram to a matrix factorization is well-defined.
\end{pro}
\begin{proof}
This is obvious by Lemma \ref{com-ass}.
\end{proof}
\indent
The next Proposition from (1) to (7) claims that the inner marking built by gluing diagrams can be removed.

\begin{pro}

\begin{enumerate}
\item[{\rm (1)}] There is an isomorphism in $\HMF^{gr}_{\Q [x] \otimes R,x^{n+1} - \omega}$,
where $R$ is a polynomial ring generated by boundary parameters except $x$,
$$
\displaystyle \c\acute{}\Big( \ \input{figgluing-mf3}\ \Big)_{n} \simeq \c\acute{}\Big( \input{figgluing-mf6} \Big)_{n} .
$$
\item[{\rm (2)}] There is an isomorphism in $\HMF^{gr}_{\Q [x] \otimes R,\omega - x^{n+1}}$, 
where $R$ is a polynomial ring generated by boundary parameters except $x$,
$$
\displaystyle \c\acute{}\Big( \ \input{figgluing-mf4}\ \Big)_{n} \simeq \c\acute{}\Big( \input{figgluing-mf5} \Big)_{n} .
$$
\\
\item[{\rm (3)}] There is an isomorphism in $\HMF^{gr}_{\Q ,0}$
\begin{eqnarray*}
\c\acute{}\Big( \ \input{figcircle-mf} \ \Big)_{n} &=& K(\pi_{12};x_1-x_2)_{\Q [x_1,x_2]}\Big|_{x_1 = x_2}\\
&\simeq& (\, 0\to \Q [x]/\left< x^n \right> \{ 1-n \} \to 0 \,).
\end{eqnarray*}
\item[{\rm (4)}] There is an isomorphism in $\HMF^{gr}_{\Q [x_1,x_2,x_5,x_6] ,x_1^{n+1}+x_2^{n+1}-x_5^{n+1}-x_6^{n+1}}$
$$\c\acute{}\Bigg( \ \input{figsmoothbubble1-a}\ \Bigg)_{n} 
\simeq \c\acute{}\Bigg( \hspace{0.4cm} \input{figsmoothing2sln-mf2} \hspace{0.4cm} \Bigg)_{n}\{ -1 \} 
\oplus\c\acute{}\Bigg( \hspace{0.4cm} \input{figsmoothing2sln-mf2} \hspace{0.4cm} \Bigg)_{n}\{ 1 \}.
$$
\\
\item[{\rm (5)}] There is an isomorphism in $\HMF^{gr}_{\Q [x_1,x_2],x_1^{n+1}-x_2^{n+1}}$
$$
\displaystyle \c\acute{}\Bigg( \ \input{figsmoothbubble2-mf2}\ \Bigg)_{n} \simeq \bigoplus_{i=0}^{n-2}\c\acute{}\Bigg( \hspace{0.2cm} \input{figsmoothing4sln1-mf} \hspace{0.2cm} \Bigg)_{n} \{ 2-n+2i \}\left< 1 \right> .
$$
\item[{\rm (6)}] There is an isomorphism in $\HMF^{gr}_{\Q [x_1,x_2,x_3,x_4],x_1^{n+1}-x_2^{n+1}+x_3^{n+1}-x_4^{n+1}}$
$$
\displaystyle \c\acute{}\Bigg( \ \input{figsmoothbubble3-mf}\ \Bigg)_{n} \simeq \c\acute{}\Bigg( \ \input{figsmoothing4sln-mf}\ \Bigg)_{n} \oplus \bigoplus_{i=0}^{n-3}\c\acute{}\Bigg( \hspace{0.2cm} \input{figsmoothing3sln-mf} \hspace{0.2cm} \Bigg)_{n} \{ 3-n+2i \} .
$$
\\
\item[{\rm (7)}] There is an isomorphism in $\HMF^{gr}_{\Q [x_1,x_2,x_3,x_4,x_5,x_6],x_1^{n+1}+x_2^{n+1}+x_3^{n+1}-x_4^{n+1}-x_5^{n+1}-x_6^{n+1}}$
$$\c\acute{}\Bigg( \ \input{figsmoothingreid1-mf}\ \Bigg)_{n} \oplus \c\acute{}\Bigg( \ \input{figsmoothingreid2-mf}\ \Bigg)_{n} \simeq \c\acute{}\Bigg( \ \input{figsmoothingreid3-mf}\ \Bigg)_{n} \oplus \c\acute{}\Bigg( \ \input{figsmoothingreid4-mf}\ \Bigg)_{n}.
$$
\end{enumerate}

\end{pro}

\indent \\\\
\begin{proof}
See \cite{KR1}.
\end{proof}
The matrix factorization for the single loop 
\unitlength 0.1in
\begin{picture}(2.750,1.50)(2.2500,-4.0)
%
\special{pn 8}%
\special{ar 350 350 120 120  0.0000000 6.2831853}%
\end{picture}%
 is defined as the following matrix factorization;
$$ \c\acute{}\Big(
\unitlength 0.1in
\begin{picture}(3.0,2.50)(2.00,-4.0)
%
\special{pn 8}%
\special{ar 350 350 150 150  0.0000000 6.2831853}%
\end{picture}%
\Big)=\Big(\, 0\to \Q [x]/\left<x^n\right>\{ 1-n \} \to 0 \,\Big) .$$
\indent
Since the potential of the matrix factorization $\c\acute{}(
\unitlength 0.1in
\begin{picture}(2.750,1.50)(2.2500,-4.0)
%
\special{pn 8}%
\special{ar 350 350 120 120  0.0000000 6.2831853}%
\end{picture}%
)_n$ is $0$, we can take the homology group $\H$ of this matrix factorization. 
\begin{cor}
The Euler characteristic of the homology $\H(\c\acute{}(
\unitlength 0.1in
\begin{picture}(2.750,1.50)(2.2500,-4.0)
%
\special{pn 8}%
\special{ar 350 350 120 120  0.0000000 6.2831853}%
\end{picture}%
)_n)$ equals the value of $\mathfrak{sl}_n$ quantum link invariant for the single loop;
$$
\overline{\chi} \Big(\H(\c\acute{}\Big(
\unitlength 0.1in
\begin{picture}(3.0,2.50)(2.00,-4.0)
%
\special{pn 8}%
\special{ar 350 350 150 150  0.0000000 6.2831853}%
\end{picture}%
\Big)_n)\Big)=[n].
$$
\end{cor}
\begin{proof}
See \cite{KR1}.
\end{proof}
The map $\c\acute{}$ from a planar diagram to a matrix factorization extends
the map $\c$ from the projection diagram associated with a link to a complex of matrix factorizations.
$$
\xymatrix{
\c\Big( \input{figplus} \Big)  := \Big( \ar[r] 
&*^{\stackrel{-1 \atop {}}{\c\acute{}\Big( \input{figsmoothing2sln1} \Big)\{ n\} \left< 1 \right>}}\ar[r] 
&*^{\stackrel{0 \atop {}}{\c\acute{}\Big( \input{figsmoothing1sln} \Big)\{ n-1 \} \left< 1 \right>}} \ar[r] 
&*^{0}\ar[r] & \Big)  \\
\c\Big( \input{figminus} \Big)  := \Big(\ar[r] 
&*^{0}\ar[r] 
&*^{\stackrel{0 \atop {}}{ \c\acute{}\Big( \input{figsmoothing1sln} \Big)\{ -n+1 \} \left< 1 \right>}}\ar[r] 
&*^{\stackrel{1 \atop {}}{ \c\acute{}\Big( \input{figsmoothing2sln1} \Big)\{ -n \} \left< 1 \right>}}\ar[r] &\Big) 
}
$$
\begin{thm}\label{KRhom}{\rm [Khovanov-Rozansky,\cite{KR1}]}
The map $\c$ is an oriented link invariant in the homotopy category of a complex of matrix factorizations.
\end{thm}
\begin{proof}
See \cite{KR1}. 
\end{proof}
%
%
%
%
\subsection{Definition of matrix factorization for trivalent diagrams and double line}\label{sec3.3}
\indent 
We extend the map $\c\acute{}$ for a planar diagram with two kinds of the boundaries, 
which consist of end points of a single line and a double line.
We define matrix factorizations for the essential planar diagrams \input{figdline2}, \input{fig3valent-in} and \input{fig3valent-out} as follows.\\
\indent
After this, we assume that $n \geq 3$. We assign the parameter $x_i$ on the end point of a single line 
and the pair of parameters $(y_j,z_j)$ on the end point of a double line. 
And each $\Z$-grading of $x_i$ and $y_j$ also equals $2$, but the $\Z$-grading of $z_i$ equals $4$;
$$
\deg(x_i)=2{\rm , \hspace{2mm}}\deg(y_i)=2 {\rm \hspace{2mm} and \hspace{2mm}} \deg(z_i)=4.
$$
Define a matrix factorization for the first diagram, double line, by
$$
\c\acute{}\Bigg(\input{figdline-mf}\Bigg)_n=K \left(\left(
\begin{array}{c}
\frac{f(y_1,z_1) - f(y_2,z_1)}{y_1 - y_2}\\
\frac{f(y_2,z_1) - f(y_2,z_2)}{z_1 - z_2}
\end{array}
\right);
\left(
\begin{array}{c}
y_1 - y_2\\
z_1 - z_2
\end{array}
\right)
\right)_{\Q [y_1,y_2,z_1,z_2]}
$$
 in $\MF^{gr}_{\Q [y_1,y_2,z_1,z_2],f(y_1,z_1)-f(y_2,z_2)}$.\\
Define a matrix factorization for the second trivalent diagram by
$$
\c\acute{}\Bigg(\input{fig3valent-in-mf}\Bigg)_n=K\left(\left(
\begin{array}{c}
\frac{f(y_3,z_3) - f(x_1 + x_2 ,z_3)}{y_3 - x_1 - x_2}\\
\frac{f(x_1+x_2,z_3) - f(x_1 + x_2,x_1 x_2)}{z_3 - x_1 x_2}
\end{array}
\right);
\left(
\begin{array}{c}
y_3 - x_1 - x_2\\
z_3 - x_1 x_2
\end{array}
\right)
\right)_{\Q [x_1,x_2,y_3,z_3]}
$$
 in $\MF^{gr}_{\Q [x_1,x_2,y_3,z_3],f(y_3,z_3)-x_1^{n+1}-x_2^{n+1}}$.\\
Define a matrix factorization for the third trivalent diagram by
$$
\c\acute{}\Bigg(\input{fig3valent-out-mf}\Bigg)_n=K\left(\left(
\begin{array}{c}
\frac{f(x_1+x_2,x_1 x_2) - f(y_3,x_1 x_2)}{x_1 + x_2 - y_3}\\
\frac{f(y_3,x_1 x_2) - f(y_3,z_3)}{x_1 x_2 - z_3}
\end{array}
\right);
\left(
\begin{array}{c}
x_1 + x_2 - y_3\\
x_1 x_2 - z_3
\end{array}
\right)
\right)_{\Q [x_1,x_2,y_3,z_3]} \{ -1 \}
$$
 in $\MF^{gr}_{\Q [x_1,x_2,y_3,z_3],x_1^{n+1}+x_2^{n+1}-f(y_3,z_3)}$.\\
\indent
The basic potential of a double line is the polynomial $f(y,z)$ obtained by expanding the power sum with the elementary symmetric polynomials.
This polynomial $f(y,z)$ is a non-homogeneous polynomial, but it has the homogeneous $\Z$-grading $2n+2$.\\
\indent
Now, we consider two planar diagrams, which can match at end points of oriented double lines with keeping the orientation.
For the matrix factorizations 
$
\c\acute{}\Big(\input{figdgluing-mf1}\Big)_n
\in \ob (\MF^{gr}_{R,\omega + f(y,z)})
$
 and 
$
\c\acute{}\Big(\input{figdgluing-mf2}\Big)_n
\in \ob( \MF^{gr}_{R\acute{},\omega\acute{} - f(y\acute{},z\acute{})}),
$
the matrix factorization $\c\acute{}\Big(\input{figdgluing-mf3}\Big)_n$ is defined by the tensor product 
\begin{center}
$\c\acute{}\Big(\input{figdgluing-mf1}\Big)_n \boxtimes \c\acute{}\Big(\input{figdgluing-mf2}\Big)_n \Big|_{(y\acute{},z\acute{})=(y,z)}$ 
$\in \ob( \MF^{gr}_{R \otimes R\acute{},\omega + \omega\acute{}})$.
\end{center}
\indent
For a polynomial $f \in \Q [x_1,x_2,\cdots ,x_k]$, the Jacobi algebra $J_f$ is defined 
as the quotient ring
$$
\Q [x_1,x_2,\cdots ,x_k]/\left<\frac{\partial f}{\partial x_1},\frac{\partial f}{\partial x_2},\cdots ,\frac{\partial f}{\partial x_k}\right> ,
$$
where $\left<\frac{\partial f}{\partial x_1},\frac{\partial f}{\partial x_2},\cdots ,\frac{\partial f}{\partial x_k}\right>$ is 
the ideal of $\Q [x_1,x_2,\cdots ,x_k]$ generated by $\frac{\partial f}{\partial x_1},\frac{\partial f}{\partial x_2},\cdots ,\frac{\partial f}{\partial x_k}$.
\begin{pro}\label{tetra}
\begin{enumerate}
\item[{\rm (1)}] There is an isomorphism in $\HMF^{gr}_{\Q [y,z]\otimes R,f(y,z)-\omega}$, 
where R is a polynomial ring generated by boundary parameters except $y$ and $z$,
$$
\c\acute{}\Big({\input{figdgluing-mf4}}\Big)_n \simeq \c\acute{}\Big({\input{figdgluing-mf1}}\Big)_n.
$$
\item[{\rm (2)}] There is an isomorphism in $\HMF^{gr}_{\Q [y,z] \otimes R,\omega-f(y,z)}$, 
where R is a polynomial ring generated by boundary parameters except $y$ and $z$,
$$
\c\acute{}\Big({\input{figdgluing-mf5}}\Big)_n \simeq \c\acute{}\Big({\input{figdgluing-mf2}}\Big)_n.
$$
\\
\item[{\rm (3)}] There is an isomorphism in $\HMF^{gr}_{\Q ,0}$
\begin{eqnarray*}
\c\acute{}\Big({\input{figdcircle-mf}}\Big)_n &=& K \left(\left(
\begin{array}{c}
\frac{f(y_1,z_1) - f(y_2,z_1)}{y_1 - y_2}\\
\frac{f(y_2,z_1) - f(y_2,z_2)}{z_1 - z_2}
\end{array}
\right);
\left(
\begin{array}{c}
y_1 - y_2\\
z_1 - z_2
\end{array}
\right)
\right)_{\Q [y_1,y_2,z_1,z_2]}\Big|_{(y_1,z_1)=(y_2,z_2)}\\
&\simeq& 
\Big(
\xymatrix{
*{J_{f(y,z)}\{ 4-2n\}}\ar[rr]&&0\ar[rr]&&*^{J_{f(y,z)}\{ 4-2n\}} 
}
\Big)
.
\end{eqnarray*}
\item[{\rm (4)}]There is an isomorphism in $\HMF^{gr}_{\Q [y_1,y_2,z_1,z_2],f(y_1,z_1)-f(y_2,z_2)}$
$$
\c\acute{}\Bigg({\input{figdgluing1}}\Bigg)_n \simeq \c\acute{}\Bigg({\input{figdgluing2}}\Bigg)_n\{ -1\} \oplus \c\acute{}\Bigg({\input{figdgluing2}}\Bigg)_n\{ 1\}.
$$
\\
\item[{\rm (5)}]There is an isomorphism in $\HMF^{gr}_{\Q [x_1,x_2,x_3,x_4],x_1^{n+1}+x_2^{n+1}-x_3^{n+1}-x_4^{n+1}}$
$$
\c\acute{}\Bigg({\input{figdgluing3}}\Bigg)_n \simeq \c\acute{}\Bigg({\input{figdgluing4}}\Bigg)_n .
$$
\end{enumerate}
\vspace{1cm}
\end{pro}
\begin{proof}
(1)
The matrix factorization $\c\acute{}\Bigg({\input{figdgluing-mf6}}\Bigg)_n$ is described by
$$
(\xymatrix{M_0 \ar[r] & M_1 \ar[r] & M_0}) ,
$$
where $M_i$ is an $R \otimes \Q [x\acute{} , y\acute{}]$-module. Then, we have
$$
\c\acute{}\Bigg({\input{figdgluing-mf4}}\Bigg)_n 
= (\xymatrix{M_0 \ar[r] & M_1 \ar[r] & M_0}) \boxtimes 
K\left( \left(
\begin{array}{c}
\frac{f(y,z)-f(y\acute{},z)}{y-y\acute{}}\\
\frac{f(y\acute{},z)-f(y\acute{},z\acute{})}{z-z\acute{}}
\end{array}
\right) ;
\left(
\begin{array}{c}
y-y\acute{}\\
z-z\acute{}
\end{array}
\right)
\right)_{\Q [y , y\acute{} , z , z\acute{}]}
$$
Since the potential of this matrix factorization does not include the parameters $y\acute{}$ and $z\acute{}$,
we can apply Theorem \ref{exclude} to $y-y\acute{}$ and $z-z\acute{}$.
Hence, we obtain the following isomorphic
\begin{eqnarray*}
&&(\xymatrix{M_0 \ar[r] & M_1 \ar[r] & M_0}) \boxtimes 
K\left( \left(
\begin{array}{c}
\frac{f(y,z)-f(y\acute{},z)}{y-y\acute{}}\\
\frac{f(y\acute{},z)-f(y\acute{},z\acute{})}{z-z\acute{}}
\end{array}
\right) ;
\left(
\begin{array}{c}
y-y\acute{}\\
z-z\acute{}
\end{array}
\right)
\right)_{\Q [y , y\acute{} , z , z\acute{}]}\\
&&\hspace{1cm}\simeq
(\xymatrix{*^{M_0/\left< y-y\acute{},z-z\acute{}\right> M_0}  \ar[r] &*^{M_1/\left< y-y\acute{},z-z\acute{}\right> M_1} \ar[r] &*^{M_0/\left< y-y\acute{},z-z\acute{}\right> M_0}})\\
&&\hspace{1cm}\simeq
\c\acute{}\Big({\input{figdgluing-mf1}}\Big)_n .
\end{eqnarray*}
(2)This proof is similar to (1) \\
(3) 
The polynomial $f(y_1,z_1)$ is explicitly described as
$$
f(y_1,z_1)= y_1^{n+1}+(n+1)\sum_{1 \leq 2i \leq n+1} \frac{(-1)^i}{i}
\left(
\begin{array}{c}
n-i\\
i-1
\end{array}
\right)
y_1^{n+1-2i} z_1^i
.$$
Then, we have
\begin{eqnarray*}
\frac{\partial f(y_1,z_1)}{\partial y_1}&=& (n+1)y_1^{n}+(n+1)\sum_{1 \leq 2i \leq n} \frac{(-1)^i(n+1-2i)}{i}\left(
\begin{array}{c}
n-i\\
i-1
\end{array}
\right)
y_1^{n-2i} z_1^i ,\\
\frac{\partial f(y_1,z_1)}{\partial z_1}&=&(n+1)\sum_{1 \leq 2i \leq n+1} (-1)^i
\left(
\begin{array}{c}
n-i\\
i-1
\end{array}
\right)
y_1^{n+1-2i} z_1^{i-1} .
\end{eqnarray*}
\indent
In the case that $n$ is even:\\
The sequence ($\frac{\partial f(y_1,z_1)}{\partial y_1},\frac{\partial f(y_1,z_1)}{\partial z_1}$) can be described as
$$
(\, (n+1)y_1^{n}+(-1)^{\frac{n}{2}}(n+1)z_1^{\frac{n}{2}} + p(y_1,z_1),-(n+1)y_1^{n-1} + q(y_1,z_1)\, )
,$$
where the polynomial degree $p(y_1,z_1)$ for $y_1$ satisfies $$\pdeg_{y_1} (p(y_1,z_1)) < n-1$$ 
and the polynomial degree for $z_1$ satisfies $$\pdeg_{z_1} (p(y_1,z_1)) < \frac{n}{2}$$ 
and the polynomial $q(y_1,z_1)$ satisfies that $$\pdeg_{y_1} (q(y_1,z_1)) < n-2$$ and $$\pdeg_{z_1} (q(y_1,z_1)) < \frac{n}{2}.$$\\
It is easy to find that this sequence is a regular sequence. By Lemma \ref{identity} and Lemma \ref{functor}, we have
\begin{eqnarray*}
\c\acute{}\Big({\input{figdcircle-mf}}\Big)_n 
&\simeq&
K \left(
\left(
\begin{array}{c}
0\\
0
\end{array}
\right);
\left(
\begin{array}{c}
-(n+1)y_1^{n}-(-1)^{\frac{n}{2}}(n+1)z_1^{\frac{n}{2}} - p(y_1,z_1)\\
(n+1)y_1^{n-1} - q(y_1,z_1)
\end{array}
\right)
\right)_{\Q [y_1,z_1]} \{ 1-n\}\{ 3-n\} 
\\
&\simeq&
\begin{array}{c}
\Big(\xymatrix{\Q [y_1,z_1]\ar[r]&0\ar[r]&\Q [y_1,z_1]}\Big)\\
\boxtimes \\
K \left(
\left(
\begin{array}{c}
0\\
0
\end{array}
\right);
\left(
\begin{array}{c}
-(n+1)y_1^{n}-(-1)^{\frac{n}{2}}(n+1)z_1^{\frac{n}{2}} - p(y_1,z_1)\\
(n+1)y_1^{n-1} - q(y_1,z_1)
\end{array}
\right)
\right)_{\Q [y_1,z_1]} 
\end{array}\{ 4-2n\} .
\end{eqnarray*}
We apply Theorem \ref{exclude} to the polynomial $(n+1)y_1^{n-1} - q(y_1,z_1)$ for $y_1$. Then, we obtain
$$
\begin{array}{c}
(\xymatrix{\Q [y_1,z_1]/\left< (n+1)y_1^{n-1} - q(y_1,z_1)\right> \ar[r]&0\ar[r]&\Q [y_1,z_1]/\left< (n+1)y_1^{n-1} - q(y_1,z_1)\right> }\\
 \boxtimes \\
K \left(
0;-(-1)^{\frac{n}{2}}(n+1)z_1^{\frac{n}{2}} - y_1 q(y_1,z_1) - p(y_1,z_1)
\right)_{\Q [y_1,z_1]/\left< (n+1)y_1^{n-1} - q(y_1,z_1)\right> }
\end{array}\{ 4-2n\} .
$$
Since we have $\Q [y_1,z_1]/\left< (n+1)y_1^{n-1} - q(y_1,z_1)\right> \simeq \Q [y_1]/\left< y_1^{n-1}\right>[z_1]$ , $\pdeg_{z_1} (y_1 q(y_1,z_1)+ p(y_1,z_1))<\frac{n}{2}$ and 
the sequence $((n+1)y_1^{n}+(-1)^{\frac{n}{2}}(n+1)z_1^{\frac{n}{2}} + p(y_1,z_1),-(n+1)y_1^{n-1} + q(y_1,z_1))$ is regular, 
we can once apply Theorem \ref{exclude} to the polynomial $-(-1)^{\frac{n}{2}}(n+1)z_1^{\frac{n}{2}} - y_1 q(y_1,z_1) - p(y_1,z_1)$ for $z_1$. 
Hence, we obtain
$$
\c\acute{}\Big( \input{figdcircle-mf} \Big)_n \simeq \Big(
\xymatrix{
*^{J_{f(y_1,z_1)}\{ 4-2n\}}
\ar[rr]&&
0
\ar[rr]&&
*^{J_{f(y_1,z_1)}\{ 4-2n\}} 
}
\Big)
.
$$\\
\indent
In the case that $n$ is odd:\\
The sequence ($\frac{\partial f(y_1,z_1)}{\partial y_1},\frac{\partial f(y_1,z_1)}{\partial z_1}$) can be described as
$$
((n+1)y_1^{n}+ p(y_1,z_1), -(n+1)y_1^{n-1}+(-1)^{\frac{n+1}{2}}(n+1)z_1^{\frac{n-1}{2}} + q(y_1,z_1))
,$$
where the polynomial degree $p(y_1,z_1)$ for $y_1$ satisfies $$\pdeg_{y_1} (p(y_1,z_1)) < n-1$$ 
and the polynomial degree for $z_1$ satisfies $$\pdeg_{z_1} (p(y_1,z_1)) < \frac{n+1}{2}$$ 
and the polynomial $q(y_1,z_1)$ satisfies that $$\pdeg_{y_1} (q(y_1,z_1)) < n-2$$ and $$\pdeg_{z_1} (q(y_1,z_1)) < \frac{n-1}{2}.$$
By Lemma \ref{identity} and Lemma \ref{functor}, we have
$$
\c\acute{}\Big({\input{figdcircle-mf}}\Big)_n 
\simeq
\begin{array}{c}
\Big(\xymatrix{\Q [y_1,z_1]\ar[r]&0\ar[r]&\Q [y_1,z_1]}\Big)\\
\boxtimes \\
K \left(
\left(
\begin{array}{c}
0\\
0
\end{array}
\right);
\left(
\begin{array}{c}
-(n+1)y_1^{n} - p(y_1,z_1)\\
(n+1)y_1^{n-1} - (-1)^{\frac{n+1}{2}}(n+1)z_1^{\frac{n-1}{2}} - q(y_1,z_1)
\end{array}
\right)
\right)_{\Q [y_1,z_1]} \{ 1-n\}\{ 3-n\} .
\end{array}
$$
Now, we have that $\pdeg_{y_1}(p(y_1,z_1))< n-1$, $\pdeg_{z_1}(q(y_1,z_1)-(n+1)y_1^{n-1}) < \frac{n-1}{2}$ 
and the sequence ($\frac{\partial f(y_1,z_1)}{\partial y_1},\frac{\partial f(y_1,z_1)}{\partial z_1}$) is regular.
Thus, after we apply Theorem \ref{exclude} to $\,-(n+1)y_1^{n} - p(y_1,z_1)\,$ for $y_1$, 
we can once apply Theorem \ref{exclude} to \\$(n+1)y_1^{n-1} - (-1)^{\frac{n+1}{2}}(n+1)z_1^{\frac{n-1}{2}} - q(y_1,z_1)$ for $z_1$. 
We also obtain 
$$
\c\acute{}\Big( \input{figdcircle-mf} \Big)_n \simeq \Big(
\xymatrix{
*^{J_{f(y_1,z_1)}\{ 4-2n\}}
\ar[rr]&&
0
\ar[rr]&&
*^{J_{f(y_1,z_1)}\{ 4-2n\}} 
}
\Big)
.
$$
(4)By definition, we have
$$
\c\acute{}\Bigg({\input{figdgluing1}}\Bigg)_n = 
K\left( \left(
\begin{array}{c}
\frac{f(y_1,z_1)-f(x_3 + x_4,z_1)}{y_1 - x_3 - x_4}\\
\frac{f(x_3 + x_4,z_1)-f(x_3 + x_4,x_3 x_4)}{z_1-x_3 x_4}\\
\frac{f(x_3+x_4,x_3 x_4)-f(y_2,x_3 x_4)}{x_3+x_4-y_2}\\
\frac{f(y_2,x_3 x_4)-f(y_2,z_2)}{x_3 x_4-z_2}
\end{array}
\right) ;
\left(
\begin{array}{c}
y_1 - x_3 - x_4\\
z_1-x_3 x_4\\
x_3+x_4-y_2\\
x_3 x_4-z_2
\end{array}
\right)
\right)_{\Q [y_1,z_1,y_2,z_2,x_3,x_4]} \{ -1\} .
$$
Since the potential of this matrix factorization does not include the parameter $x_3$ and $x_4$,
 we can apply Theorem \ref{exclude} to the polynomial $y_1 - x_3 - x_4$ for $x_3$. 
Thus, this matrix factorization is equivalent to the following matrix factorization
$$
K\left( \left(
\begin{array}{c}
\frac{f(y_1,z_1)-f(y_1,-x_4^2 +y_1 x_4)}{z_1+x_4^2 - y_1 x_4}\\
\frac{f(y_1,-x_4^2 +y_1 x_4)-f(y_2,-x_4^2 +y_1 x_4)}{y_1-y_2}\\
\frac{f(y_2,-x_4^2 +y_1 x_4)-f(y_2,z_2)}{-x_4^2 +y_1 x_4-z_2}
\end{array}
\right) ;
\left(
\begin{array}{c}
z_1+x_4^2 - y_1 x_4\\
y_1-y_2\\
-x_4^2 +y_1 x_4-z_2
\end{array}
\right)
\right)_{\Q [y_1,z_1,y_2,z_2,x_3,x_4]/\left< y_1 - x_3 - x_4 \right> } \{ -1\} 
$$
 in $\HMF^{gr}_{\Q [y_1,z_1,y_2,z_2,x_4],f(y_1,z_1)-f(y_2,z_2)}$. 
Since the quotient ring $\Q [y_1,z_1,y_2,z_2,x_3,x_4]/\left< y_1 - x_3 - x_4 \right> $ is isomorphic to $$\Q [y_1,z_1,y_2,z_2,x_4],$$
we once apply Theorem \ref{exclude} to the polynomial $z_1+x_4^2 - y_1 x_4$ for $x_4$. 
Then, this matrix factorization is equivalent to the following matrix factorization 
$$
K\left( \left(
\begin{array}{c}
\frac{f(y_1,z_1)-f(y_2,z_1)}{y_1-y_2}\\
\frac{f(y_2,z_1)-f(y_2,z_2)}{z_1-z_2}
\end{array}
\right) ;
\left(
\begin{array}{c}
y_1-y_2\\
z_1-z_2
\end{array}
\right)
\right)_{\Q [y_1,z_1,y_2,z_2,x_4]/\left< z_1+x_4^2 - y_1 x_4\right> } \{ -1\} 
$$
in $\HMF^{gr}_{\Q [y_1,z_1,y_2,z_2],f(y_1,z_1)-f(y_2,z_2)}$. 
The quotient ring $\Q [y_1,z_1,y_2,z_2,x_4]/\left< z_1+x_4^2 - y_1 x_4 \right> \{ -1\}$ is isomorphic to 
$$\Q [y_1,z_1,y_2,z_2]\{ -1 \} \oplus \Q [y_1,z_1,y_2,z_2]\, x_4 \{ -1\}.$$
Furthermore, this ring is equivalent to $$\Q [y_1,z_1,y_2,z_2]\{ -1 \} \oplus \Q [y_1,z_1,y_2,z_2] \{ 1\}$$ as a $\Z$-graded $\Q [y_1,z_1,y_2,z_2]$-module.
Hence, we obtain the following matrix factorization
$$
K\left( \left(
\begin{array}{c}
\frac{f(y_1,z_1)-f(y_2,z_1)}{y_1-y_2}\\
\frac{f(y_2,z_1)-f(y_2,z_2)}{z_1-z_2}
\end{array}
\right) ;
\left(
\begin{array}{c}
y_1-y_2\\
z_1-z_2
\end{array}
\right)
\right)_{\Q [y_1,z_1,y_2,z_2]} \{ -1\}
\oplus
K\left( \left(
\begin{array}{c}
\frac{f(y_1,z_1)-f(y_2,z_1)}{y_1-y_2}\\
\frac{f(y_2,z_1)-f(y_2,z_2)}{z_1-z_2}
\end{array}
\right) ;
\left(
\begin{array}{c}
y_1-y_2\\
z_1-z_2
\end{array}
\right)
\right)_{\Q [y_1,z_1,y_2,z_2]} \{ 1\} .
$$
This is the right-hand side of the equivalent of Proposition ($4$).\\
(5)By definition, we have
$$
\c\acute{}\Bigg({\input{figdgluing3}}\Bigg)_n = 
K\left( \left(
\begin{array}{c}
\frac{f(y_5,z_5)-f(x_3 + x_4,z_5)}{y_5 - x_3 - x_4}\\
\frac{f(x_3 + x_4,z_5)-f(x_3 + x_4,x_3 x_4)}{z_5-x_3 x_4}\\
\frac{f(x_1+x_2,x_1 x_2)-f(y_5,x_1 x_2)}{x_1+x_2-y_5}\\
\frac{f(y_5,x_1 x_2)-f(y_5,z_5)}{x_1 x_2-z_5}
\end{array}
\right) ;
\left(
\begin{array}{c}
y_5 - x_3 - x_4\\
z_5-x_3 x_4\\
x_1+x_2-y_5\\
x_1 x_2-z_5
\end{array}
\right)
\right)_{\Q [x_1,x_2,x_3,x_4,y_5,z_5]} \{ -1\} .
$$
Since the potential of this matrix factorization does not include the parameter $y_5$ and $z_5$, 
we apply Theorem \ref{exclude} to the polynomials $y_5 - x_3 - x_4$ and $z_5-x_3 x_4$.
Then, it is easy to find that we obtain the right-hand side of the equivalent of Proposition (5).
\end{proof}
A matrix factorization for the double loop 
\unitlength 0.1in
\begin{picture}(3.30,1.500)(1.85,-3.75)
%
\special{pn 8}%
\special{ar 350 350 120 120  0.0000000 6.2831853}%
%
\special{pn 8}%
\special{ar 350 350 100 100  0.0000000 6.2831853}%
\end{picture}%
 is defined as the above matrix factorization;
$$ 
\c\acute{}\Big(
\unitlength 0.1in
\begin{picture}(3.30,2.500)(1.85,-3.75)
%
\special{pn 8}%
\special{ar 350 350 150 150  0.0000000 6.2831853}%
%
\special{pn 8}%
\special{ar 350 350 130 130  0.0000000 6.2831853}%
\end{picture}%
\Big)_n = \Big( \xymatrix{*{J_{f(y,z)}\{ 4-2n\}}\ar[rr]&&0\ar[rr]&&*^{J_{f(y,z)}\{ 4-2n\}} }\Big) . 
$$
\indent
Since the potential of the matrix factorization $\c\acute{}(
\unitlength 0.1in
\begin{picture}(3.30,1.500)(1.85,-3.75)
%
\special{pn 8}%
\special{ar 350 350 120 120  0.0000000 6.2831853}%
%
\special{pn 8}%
\special{ar 350 350 100 100  0.0000000 6.2831853}%
\end{picture}%
)_n $ is $0$, we can take the homology group $\H$ of this matrix factorization. 
\begin{cor}\label{dloop}
The Euler characteristic of the homology $\H (\c\acute{}(
\unitlength 0.1in
\begin{picture}(3.30,1.500)(1.85,-3.75)
%
\special{pn 8}%
\special{ar 350 350 120 120  0.0000000 6.2831853}%
%
\special{pn 8}%
\special{ar 350 350 100 100  0.0000000 6.2831853}%
\end{picture}%
)_n)$ equals the value of $\mathfrak{sl}_n$ quantum link invariant for the double loop;
$$
\overline{\chi}\Big( \H (\c\acute{}\Big(
\unitlength 0.1in
\begin{picture}(3.30,2.500)(1.85,-3.75)
%
\special{pn 8}%
\special{ar 350 350 150 150  0.0000000 6.2831853}%
%
\special{pn 8}%
\special{ar 350 350 130 130  0.0000000 6.2831853}%
\end{picture}%
\Big)_n) \Big) =\frac{[n][n-1]}{[2]}.
$$
\end{cor}
\begin{proof}
By proof of the above Proposition \ref{tetra} (3), we have
$$
J_{f(y,z)}\simeq \left\{
\begin{array}{cc}
\Q [y,z]/\left< y^{n-1},z^{\frac{n}{2}}\right> & n:{\rm even}\\
\\
\Q [y,z]/\left< y^{n},z^{\frac{n-1}{2}}\right> & n:{\rm odd}
\end{array}
\right.
$$
(i) $n$ : even\\
Since $\deg(y_1)=2$ and $\deg(z_1)=4$, we have
\begin{eqnarray*}
\overline{\chi}\Big( \H ( \c\acute{}\Big(
\unitlength 0.1in
\begin{picture}(3.30,2.500)(1.85,-3.75)
%
\special{pn 8}%
\special{ar 350 350 150 150  0.0000000 6.2831853}%
%
\special{pn 8}%
\special{ar 350 350 130 130  0.0000000 6.2831853}%
\end{picture}%
\Big)_n ) \Big) &=&\overline{\chi}\left( J_{f(y,z)}\{ 4-2n\} \right)\\
&=& \sum_{i=0}^{n-2}\sum_{j=0}^{\frac{n-2}{2}} q^{4-2n+2i+4j}\\
&=& \sum_{i=0}^{n-2} q^{2-n+2i} \sum_{j=0}^{\frac{n-2}{2}} q^{2-n+4j}\\
&=& (q^{2-n}+q^{4-n}+q^{6-n}+\cdots+q^{n-4}+q^{n-2})(q^{2-n}+q^{6-n}+q^{10-n}+\cdots+q^{n-6}+q^{n-2})\\
&=& \frac{[n][n-1]}{[2]}.
\end{eqnarray*}
(ii) $n$ : odd\\
We have
\begin{eqnarray*}
\overline{\chi}\Big( \H ( \c\acute{}\Big(
\unitlength 0.1in
\begin{picture}(3.30,2.500)(1.85,-3.75)
%
\special{pn 8}%
\special{ar 350 350 150 150  0.0000000 6.2831853}%
%
\special{pn 8}%
\special{ar 350 350 130 130  0.0000000 6.2831853}%
\end{picture}%
\Big)_n ) \Big) &=&\overline{\chi}\left( J_{f(y,z)}\{ 4-2n\} \right)\\
&=& \sum_{i=0}^{n-1}\sum_{j=0}^{\frac{n-3}{2}} q^{4-2n+2i+4j}\\
&=& \sum_{i=0}^{n-1} q^{1-n+2i} \sum_{j=0}^{\frac{n-3}{2}} q^{3-n+4j}\\
&=& (q^{1-n}+q^{3-n}+q^{5-n}+\cdots+q^{n-3}+q^{n-1})(q^{3-n}+q^{7-n}+q^{11-n}+\cdots+q^{n-7}+q^{n-3})\\
&=& \frac{[n][n-1]}{[2]}.
\end{eqnarray*}
\end{proof}
\end{document}